\documentclass[11pt,citesort,rotate,psfig]{article}  %  ?{article}
\usepackage{epsfig,amsfonts}
\catcode`\"=13      % aktivira znak "
\def"#1{\v{#1}}      % sedaj je "c lahko \v{c}

\def\be{\begin{equation}}
\def\ee{\end{equation}}
\def\bea{\begin{eqnarray}}
\def\eea{\end{eqnarray}}

\def\ppp#1 {#1^{\prime \prime \prime}}

\def\hfl#1#2{\smash{\mathop{\hbox to 12 mm{\rightarrowfill}} \limits^{\scriptstyle #1}_{\scriptstyle #2}}}

\def\binom#1#2{\left ( {{#1} \atop {#2}} \right )}

\newcommand {\qed} {\null \hfill \rule {2.5mm}{2.5mm}}

\newcommand {\N} {\ensuremath{\mathbb{N}}}

\author{Tomislav Do\v{s}li\'c }
\title{Log-balanced combinatorial sequences}
\date{\today}

\begin{document}
\maketitle
Texas A\&M University at Galveston, Galveston, Texas, 77551, U.S.A. and \\
Faculty of Agriculture, University of Zagreb, Sveto\v simunska c. 25, 10000 Zagreb, CROATIA
\\ \\
\\

Proposed running title: Log-balanced combinatorial sequences

\newpage
{\bf \Large Abstract}

In this paper we consider log-convex sequences that satisfy an additional
constraint imposed on their rate of growth. We call such sequences 
log-balanced. It is shown that all such sequences satisfy a pair of double
inequalities. Sufficient conditions for log-balancedness are given for the
case when the
sequence satisfies a two- (or more-) term linear recurrence. It is shown
that many combinatorially interesting sequences belong to this class, and, as
a consequence, that the above mentioned double inequalities are valid for all of them.
\\ \\



\newpage
%------------------------------------------------------------
\section{Introduction}
%------------------------------------------------------------

One of the most common tasks in combinatorics is to find explicitly the size 
of a certain finite set, depending on an integer parameter $n$ and defined in
an intricate way. Then the next question usually asks how the sequence of
numbers describing this size behaves for large values of $n$. Of particular
interest is logarithmic behavior of the sequence (i.e. its log-convexity
or log-concavity), since it is often instrumental in obtaining its growth
rate and asymptotic behavior. Also, log-behavior may qualify (or disqualify)
a sequence as a candidate for use in certain models. A good example is the
recent use of log-convex sequences in quantum physics for constructing
generalized coherent states associated with models having discrete non-linear
spectra (\cite{penson}).

The literature on log-behavior of combinatorial sequences is vast; we refer
the reader to the book \cite{karlin}, and also to \cite{brenti}, \cite{stanley89}
and \cite{stanley00}. 

In this article we quantitatively refine the concept of log-convexity by 
introducing and considering the class of log-balanced
combinatorial sequences and showing that the terms of such sequences satisfy
certain double inequalities. We further proceed by deriving sufficient 
conditions for a (combinatorial) sequence given by a two-term linear 
homogeneous recurrence to be log-convex and log-balanced. It is also indicated 
how to extend this approach to longer recurrences and how to treat the case of
nonhomogeneous recurrences. Finally, we demonstrate that the class of 
log-balanced sequences is rich enough to include many cases of special 
combinatorial interest. As a consequence, we obtain new pairs of inequalities 
for many classical sequences.

\section{Log-balanced sequences}

A sequence $\left ( a_n \right )_{n \geq 0}$ of positive real numbers is
{\bf log-convex} if $a_n^2 \leq a_{n-1} a_{n+1}$ for all $n \geq 1$. If the
opposite inequality, $a_n^2 \geq a_{n-1} a_{n+1}$ is valid for all $n \geq 1$,
we say that the sequence $\left ( a_n \right )_{n \geq 0}$ is {\bf log-concave}.
In case of equality, $a_n^2 = a_{n-1} a_{n+1}$, $n \geq 1$, we call the
sequence $\left ( a_n \right )_{n \geq 0}$ geometric or {\bf log-straight}.
Another type of logarithmic behavior is that of the Fibonacci sequence, where
direction of the inequality depends on the parity of $n$. We call such
sequences {\bf log-Fibonacci}.

An alternative way of characterizing the log-behavior of a sequence is via the
sequence of quotients of its successive terms. We call the sequence
$\left ( x_n \right )_{n \geq 1}$, $x_n = \frac{a_n}{a_{n-1}}$ the 
{\bf quotient sequence} of the sequence $\left ( a_n \right )_{n \geq 0}$.
Obviously, the sequence $\left ( a_n \right )_{n \geq 0}$ is
log-convex if and only if its quotient sequence is non-decreasing. Similarly,
$\left ( a_n \right )_{n \geq 0}$ is log-concave if and only if its quotient sequence is non-increasing,
and log-Fibonacci if and only if no three successive elements of the quotient
sequence form a monotone subsequence.

In what follows, we consider log-convex sequences whose quotient sequence
does not grow too fast. We shall also assume that $a_0 = 1$, unless explicitly
stated otherwise. This restriction is not too severe, since in many combinatorially
interesting cases we put $a_0 = 1$ by convention.

A sequence $\left ( a_n \right )_{n \geq 0}$ of positive real numbers is
{\bf log-balanced} if $\left ( a_n \right )_{n \geq 0}$ is log-convex and the
sequence $\left ( \frac{a_n}{n!} \right )_{n \geq 0}$ is log-concave. In terms
of quotient sequences, this means that $x_n \leq x_{n+1} \leq \frac{n+1}{n} x_n$,
for all $n \geq 1$.

The motivation for considering such sequences comes from the recent article
\cite{asai}, where it was shown that the sequences of Bell numbers of any
order are of
this type. Since this property makes them suitable for providing important
examples in white noise theory (\cite{kuo}), it is of interest to see 
whether there are some other such sequences and to characterize them.

We start by stating in terms of log-balanced sequences the following 
observation, made in (\cite{asai}). The proof is reproduced here 
for the reader's convenience.

{\bf Proposition 1}\\ Let $\left ( a_n \right )_{n \geq 0}$ be a log-balanced
sequence. Then \\
(a) $a_n^2 \leq a_{n-1} a_{n+1} \leq \left ( 1 + \frac{1}{n} \right ) a_n^2,
\quad n \geq 1;$\\
(b) $a_n a_m \leq a_{n+m} \leq \binom{n+m}{n} a_n a_m, \quad n, m \geq 0.$

{\bf Proof}\\ The double inequality (a) is just another way of stating the
fact that the sequence $\left ( a_n \right )_{n \geq 0}$ is log-balanced.

The left inequality of (b) follows easily (by induction) from the log-convexity of 
$\left ( a_n \right )_{n \geq 0}$. To prove the right inequality, start from
$x_{n}\geq \frac{n}{n+1}x_{n+1}.$
By using this inequality repeatedly, we get
$$\frac{a_{1}}{a_{0}}\geq \frac{1}{2}\frac{a_{2}}{a_{1}}\geq \frac{1}{3}\frac{%
a_{3}}{a_{2}}\geq ...\geq \frac{1}{m+n}\frac{a_{m+n}}{a_{m+n-1}},$$ for all $%
n\geq 0$, $m\geq 1.$

Hence, for any $0\leq j\leq m-1,$ we have $$\frac{a_{j+1}}{a_{j}}\geq \frac{%
j+1}{m+n}\frac{a_{m+n}}{a_{m+n-1}}.$$
From this we get
$$\frac{a_{1}}{a_{0}}\frac{a_{2}}{a_{1}}\frac{a_{3}}{a_{2}}...\frac{a_{m}}{%
a_{m-1}}\geq \left(\frac{1}{n+1}\frac{a_{n+1}}{a_{n}}\right)\left(\frac{2}{n+2}\frac{a_{n+2}%
}{a_{n+1}}\right)...\left(\frac{m}{m+n}\frac{a_{m+n}}{a_{n}}\right).$$
After the cancellations we get $$\frac{a_{m}}{a_{0}}\geq \frac{n!m!}{(m+n)!}%
\frac{a_{m+n}}{a_{n}},$$and, taking into account the fact that $a_{0}=1,$ we
finally get
$$a_{m+n}\leq \binom{m+n}{n}a_{n}a_{m}.$$The case $m=0$ is trivially valid
for all $n\geq 0.$ \qed

\section{Sufficient conditions}

For most sequences of combinatorial interest there are no explicit, closed
form expressions for their elements. On the other hand, one can often find 
recurrences and/or generating functions for them. So, direct ways of 
establishing the log-behavior
of a given sequence (i.e. of proving inequalities of the type (a) from
Proposition 1) are only rarely at our disposal. Combinatorial proofs, which are 
the most desirable, often turn out to be rather involved and/or tricky. (A nice
survey of inductive and injective proofs of log-concavity is given in
\cite{sagan}.) Hence, it makes sense to seek analytical methods sufficiently 
robust,
easy to apply and that will work for a reasonably broad class of sequences.
Here we present one such method that works almost automatically for sequences
given by recurrence relations. We start by explaining the method for the case
of linear homogeneous recurrences of second order, and later we indicate how 
to modify this so that it can be applied also on longer and/or nonhomogeneous
recurrences.

Let $\left ( a_n \right )_{n \geq 0}$ be a sequence of positive
real numbers, given by the two-term recurrence
\be
a_n = R(n) a_{n-1} + S(n) a_{n-2}, \quad n \geq 2,
\ee
with given initial conditions $a_0$, $a_1$. The quotient sequence 
$\left ( x_n \right )_{n \geq 1}$ satisfies the nonlinear recurrence 
\be
x_n = R(n) + \frac{S(n)}{x_{n-1}}, \quad n \geq 2,
\ee
with the initial condition $x_1 = a_1 / a_0$. We assume that the sequence
$(x_n)$ is bounded by two known sequences, i.e. that there are sequences
$(m_n)$ and $(M_n)$ such that $0 < m_n \leq x_n \leq M_n$, for all $n \in \N$.
The sequences $(m_n)$ and $(M_n)$ can usually be rather easily inferred from
recurrence (2), or guessed from the initial behavior of the sequence $(x_n)$,
and then the bounding relations are verified by induction. In many cases even
the constant sequences $m_n = m$ and $M_n = M$ will be sufficiently good
lower and upper bounds for $x_n$.

As the log-convexity is of considerable interest on its own, we first establish
sufficient conditions for a sequence $(a_n)$ given by (1) to be log-convex.
We assume $R(n) \geq 0$ and treat the cases $S(n) \leq 0$ and $S(n) \geq 0$
separately. The case $S(n) \leq 0$ is simpler and we consider it first.

Assume, inductively, that $x_{n_0} \leq x_{n_0+1} \leq \ldots \leq x_n$
for some $n_0 \in \N$. Expressing $x_{n+1}$ from equation (2) and taking
into account that $S(n+1)/x_n \geq S(n+1)/x_{n-1}$, we obtain
$$x_{n+1} = R(n+1) + \frac{S(n+1)}{x_n} \geq R(n+1) + \frac{S(n+1)}{x_{n-1}}.$$
We want to prove that $x_{n+1} \geq x_n$. But this will follow if we prove
the stronger inequality in which $x_{n+1}$ is replaced by the right hand
side in the above inequality. Hence, consider the circumstance
$$R(n+1) + \frac{S(n+1)}{x_{n-1}}\geq x_n = R(n) + \frac{S(n)}{x_{n-1}},$$
or, equivalently, 
$$[R(n+1)-R(n)] x_{n-1} + S(n+1)-S(n) \geq 0.$$
By denoting $R(n+1)-R(n) = \nabla R(n)$ and $S(n+1)-S(n) = \nabla S(n)$, we get
a compact expression for the sufficient condition for the sequence $(a_n)$
to be log-convex:
\be
\nabla R(n) x_{n-1} + \nabla S(n) \geq 0, \quad n \geq n_0,
\ee
for some $n_0 \in \N$. Hence, we have established the following result:

{\bf Proposition 2}\\ Let $(a_n)_{n \geq 0}$ be a sequence of positive real 
numbers given by the two-term recurrence (1), and $(x_n)_{n \geq 1}$ its
quotient sequence, given by (2). If there is an $n_0 \in \N$ such that 
$x_{n_0} \leq x_{n_0+1}$, $R(n) \geq 0$, $S(n) \leq 0$, and
$$\nabla R(n) x_{n-1} + \nabla S(n) \geq 0,$$
for all $n \geq n_0$, then the sequence $(a_n)_{n \geq n_0}$ is log-convex. \qed

When (as is a common case) the function $R(n)$ is 
non-decreasing, the condition (3) can be further simplified without significant
loss of generality by assuming $\nabla R(n) \geq 0$ and replacing $x_{n-1}$
by $m_{n-1}$, or even by a constant $m$:
\be
\nabla R(n) m + \nabla S(n) \geq 0, \quad n \geq n_0.
\ee

The case $S(n) \geq 0$ is a bit more complicated. Again, we start from the 
inductive assumption $x_{n_0} \leq x_{n_0+1} \leq \ldots \leq x_n$ and 
want to show that $x_{n+1} \geq x_n$. By expressing both sides of this
inequality via (2), we obtain
$$R(n+1) + \frac{S(n+1)}{x_n} \geq R(n) + \frac{S(n)}{x_{n-1}}.$$
This is equivalent to
$$ x_n x_{n-1}\nabla R(n) +S(n+1) x_{n-1} -S(n) x_n \geq 0.$$
By adding the term $S(n) x_{n-1} - S(n) x_{n-1}$ to the left hand side of the
above inequality and rearranging it, we obtain
$$\nabla R(n) x_n x_{n-1} + \nabla S(n) x_{n-1} \geq S(n) (x_n - x_{n-1}).$$
Expressing the term $x_n - x_{n-1}$ via (2) now yields
$$x_{n-1} [\nabla R(n) x_n + \nabla S(n) ] \geq S(n) \left [ \nabla R(n-1)
+ \frac{S(n)}{x_{n-1}} - \frac{S(n-1)}{x_{n-2}} \right ] .$$
Now, replacing $\frac{S(n)}{x_{n-1}}$ with $\frac{S(n)}{x_{n-2}}$ in the right
hand side square brackets we get a stronger inequality which can be written as
\be
x_{n-1} x_{n-2} [ \nabla R(n) x_n + \nabla S(n) ] \geq S(n) [ \nabla R(n-1) x_{n-2} + \nabla S(n-1) ] .
\ee
Obviously, this inequality implies $x_{n+1} \geq x_n$, and it can serve as a
sufficient condition of log-convexity for the sequence $(a_n)$.

{\bf Proposition 3}\\ Let $(a_n)_{n \geq 0}$ be a sequence of positive real 
numbers given by the two-term recurrence (1), and $(x_n)_{n \geq 1}$ its
quotient sequence, given by (2). If there is an $n_0 \in \N$ such that 
$x_{n_0} \leq x_{n_0+1}$, $R(n) \geq 0$, $S(n) \geq 0$, and the inequality
$$x_{n-1} x_{n-2} [ \nabla R(n) x_n + \nabla S(n) ] \geq S(n) [ \nabla R(n-1
) x_{n-2} + \nabla S(n-1) ]$$
is valid for all $n \geq n_0$, then the sequence $(a_n)_{n \geq n_0}$
is log-convex. \qed

Again, in many combinatorially relevant cases where $\nabla R(n) \geq 0$ and
$m \leq x_n \leq M$, the sufficient condition of Proposition 3 can be simplified
to 
\be
m^2 [m \nabla R(n) + \nabla S(n) ] \geq S(n) [ M \nabla R(n-1)
+\nabla S(n-1) ].
\ee

Typically, propositions 1 and/or 2 are applied so that the respective 
inequalities are verified inductively for all $n \in \N$ greater than some
$n_0$, and the remaining cases are then checked by hand or using some 
computer algebra system.

Now we turn our attention to the inequality $x_{n+1} \leq \frac{n+1}{n} x_n$.
Again, we assume $R(n) \geq 0$ and treat the cases $S(n) \leq 0$ and $S(n)
\geq 0$ separately. Also, we assume that the log-convexity of the sequence
$(a_n)$ is already established, i.e. that the sequence $(x_n)$ is increasing.

We first consider the simpler case $R(n) \geq 0$, $S(n) \geq 0$ and find
the sufficient conditions for $x_{n+1} \leq \frac{n+1}{n} x_n$ as follows. 
From the recurrence (2)
we have $$x_{n+1} = R(n+1) + \frac{S(n+1)}{x_n}.$$
Since the sequence $\left ( x_n \right )_{n \geq 1}$ is non-decreasing, we have
$$x_{n+1} \leq R(n+1) + \frac{S(n+1)}{x_{n-1}}.$$
The condition that the right-hand side does not exceed $\frac{n+1}{n} x_n$ is
given by
$$R(n+1) + \frac{S(n+1)}{x_{n-1}} \leq \frac{n+1}{n} \left (
R(n) + \frac{S(n)}{x_{n-1}} \right ),$$
and this is equivalent to
$$n R(n+1) x_{n-1} + n S(n+1) \leq (n+1) R(n) x_{n-1} + (n+1) S(n).$$
Denoting
$$\Delta _R (n) = \left | {{R(n)} \atop {R(n+1)}} \quad {{n} \atop {n+1}} 
\right |, \quad 
\Delta _S (n) = \left | {{S(n)} \atop {S(n+1)}} \quad {{n} \atop {n+1}} \right |,$$
we get our sufficient conditions in the form 
$$ \Delta _R (n) x_{n-1} + \Delta _S (n) \geq 0.$$
Hence, we have established the following result:

{\bf Proposition 4}\\ Let $\left ( a_n \right )_{n \geq 0}$ be a log-convex 
sequence of positive real numbers given by the two-term recurrence (1).
If there is an $n_0 \in \N$ such that $x_{n_0+1} \leq \frac{n_0+1}{n_0} 
x_{n_0}$, $R(n) \geq 0$, $S(n) \geq 0$, and 
$$ \Delta _R (n) x_{n-1} + \Delta _S (n) \geq 0,$$ 
for all $n \geq n_0$, then the sequence $\left ( a_n \right )_{n \geq 0}$ 
is log-balanced.  \qed

The case $S(n) \leq 0$ is a bit more complicated. We proceed by induction on
$n$. First we check that $x_{n_0+1} \leq \frac{n_0+1}{n_0} x_{n_0}$ for some
$n_0 \in \N$, and suppose that $x_k \leq \frac{k}{k-1} x_{k-1}$ for all 
$n_0 \leq k \leq n$. Denoting $-S(n) = \tilde S(n)$, we get
$$x_{n+1} =R(n+1) - \frac{\tilde S(n+1)}{x_n}, \quad \tilde S(n+1) \geq 0.$$
From the induction hypothesis, $x_n \leq \frac{n}{n-1} x_{n-1}$, it follows
$\frac{1}{x_n} \geq \frac{n-1}{n}\frac{1}{x_{n-1}}$, and hence $-\frac{1}{x_n}
\leq - \frac{n-1}{n}\frac{1}{x_{n-1}}$. Now we have
$$x_{n+1} =R(n+1) - \frac{\tilde S(n+1)}{x_n} \leq R(n+1) - 
\frac{n-1}{n}\frac{\tilde S(n+1)}{x_{n-1}}.$$
The right hand side does not exceed $\frac{n+1}{n} x_n$ if
$$ R(n+1) -\frac{n-1}{n}\frac{\tilde S(n+1)}{x_{n-1}}
\leq \frac{n+1}{n} \left ( R(n) - \frac{\tilde S(n)}{x_{n-1}} \right ),$$
and this is, in turn, equivalent to
$$\left [ (n+1) R(n) - n R(n+1) \right ] x_{n-1} + 
(n-1)\tilde S(n+1) - (n+1) \tilde S(n) \geq 0.$$
The coefficient of $x_{n-1}$ is $\Delta _R (n)$, and the rest can be written
as
$$\left | {{n-1} \atop {n+1}} \quad {{\tilde S(n)} \atop {\tilde S(n+1)}} 
\right | = \left | {{S(n)} \atop {S(n+1)}}\quad {{n-1} \atop {n+1}} \right |.$$
Denoting the right hand side determinant by $\overline {\Delta }_S (n)$, 
we get the desired sufficient conditions:
$$ \Delta _R (n) x_{n-1} + \overline {\Delta }_S (n) \geq 0.$$
We can summarize:

{\bf Proposition 5}\\ Let $\left ( a_n \right )_{n \geq 0}$ be a log-convex 
sequence of positive real numbers given by the two-term recurrence (1) with
$R(n) \geq 0$, $S(n) \leq 0$. If there is an integer $n_0$ such that 
$x_{n_0+1} \leq \frac{n_0+1}{n_0} x_{n_0}$, and if the inequality 
$$ \Delta _R (n) x_{n-1} + \overline {\Delta }_S (n) \geq 0$$
holds for all $n \geq n_0$, then the sequence $\left ( a_n \right )_{n \geq n_0}$ is log-balanced. \qed

\section{Examples}

We now justify our introduction of log-balanced sequences by demonstrating
that the class is wide enough and that it includes many sequences of 
combinatorial relevance. As a consequence, for all our examples we establish
the validity of inequalities from Proposition 1. The left inequalities for some
of the considered sequences were established earlier (\cite{aigner}, \cite{dv}),
but the right inequalities are, with one exception (\cite{doslic}), to the best 
of our knowledge, new. For more details
on all the considered sequences, we refer the reader to the book \cite{stanleyII}
and to the references therein.

Our first example is the sequence of {\bf Motzkin numbers} (see, e.g. Ex. 6.38
of \cite{stanleyII} for its combinatorial interpretations).

{\bf Corollary 1}\\ The sequence $M_n$ of Motzkin numbers is log-balanced.

{\bf Proof}\\ The log-convexity of $M_n$ was first established
algebraically in \cite{aigner}, and a combinatorial proof appeared soon
afterwards (\cite{callan}). By our method it follows easily by starting from 
the recurrence
$$M_{n}= \frac{2n+1}{n+2}M_{n-1}+\frac{3(n-1)}{n+2}M_{n-2}, \quad n \geq 2$$
with $M_0 = M_1 = 1$. 
Here $R(n)=\frac{2n+1}{n+2} \geq 0$, $S(n)=\frac{3(n-1)}{n+2} \geq 0$. It is
easy to prove by induction on $n$ that $2 \leq M_n/M_{n-1} \leq 7/2$ for
all $n \geq 2$, and the log-convexity follows by computing $\nabla R(n)$, 
$\nabla S(n)$, $\nabla R(n-1)$, and $\nabla S(n-1)$ and then 
verifying the inequality (6). From the fact that $ \Delta _R (n) = 
\frac{2 n^2 + 4 n +3}{(n+2)(n+3)}\geq 0$, 
$ \Delta _S (n) = \frac{n^2 -n -3}{(n+2)(n+3)}\geq 0$ and $x_{n-1} \geq 0$ 
for all $n \geq 3$,
it follows that $ \Delta _R (n) x_{n-1} + \Delta _S (n) \geq 0$ for all 
$n \geq 3$. The log-balancedness of $(M_n)$
now follows from Proposition 4, after direct verification of the defining
inequality for the remaining values of $n$. \qed

Our next example is the sequence of {\bf Fine numbers}. The reader may 
consult the recent survey \cite{deutsch} for more details on Fine numbers and on
their combinatorial interpretations.

{\bf Corollary 2}\\ The sequence $B_n$ of Fine numbers is log-balanced for
$n \geq 2$.

{\bf Proof}\\ We start from the recurrence
$$B_n = \frac{7n-5}{2n+2} B_{n-1} + \frac{2n-1}{n+1} B_{n-2}, \quad n \geq 2,$$
with initial conditions $B_0 = 1$ and $B_1 = 0$. The quotient sequence, 
$x_n = B_n/B_{n-1}$, is defined for $n \geq 3$. It is easy to show, by 
induction on $n$, that $3 \leq x_n \leq 6$ for all $n \geq 3$. In fact, 
$3 \leq x_{n-1} \leq 6$ implies $3 \leq x_n \leq 6$ via the above recurrence 
for $n \geq 7$, and $x_n$ is obviously between $3$ and $6$ for $n = 3,4,5,$
and $6$. We proceed by computing $\nabla R(n) = \frac{6}{(n+1)(n+2)}$,
$\nabla S(n) = \frac{4}{(n+1)(n+2)}$, $\nabla R(n-1) = \frac{6}{n(n+1)}$, and
$\nabla S(n-1) = \frac{4}{n(n+1)}$. After plugging in these expressions we
find, condition (6) becomes $$ 10n^2 -30 n +80 \geq 0,$$
and this is true for all $n \in \N$. Hence, the sequence $(B_n)_{n \geq 2}$
is log-convex. The log-balancedness now follows by computing $ \Delta _R (n)
= \frac{7n-10}{2(n+2)}$, $ \Delta _S (n)= 2 \frac{n-1}{n+2}$, and applying
Proposition 4. \qed

The {\bf Franel numbers} of order $r$ are defined by
$$ F_n^{(r)} = \sum _{k = 0}^n \binom{n}{k} ^r.$$

{\bf Corollary 3}\\ The sequences of Franel numbers of order $3$ and $4$ are
log-balanced.

{\bf Proof}\\ It is known that Franel numbers of order $r$ satisfy a 
homogeneous linear recurrence of
order $\lfloor \frac{r+1}{2} \rfloor$ with polynomial coefficients 
(\cite{stanleyII}, p. 245-6 and p. 278). We have
$$F^{(r)}_n = R^{(r)}(n) F^{(r)}_{n-1} + S^{(r)}(n) F^{(r)}_{n-2}, \quad
r = 3, 4, \quad n \geq 2,$$
with $F_0^{(3)} = F_0^{(4)} = 1$, $F_1^{(3)} = F_1^{(4)} = 2$. Here
$$R^{(3)}(n) = \frac{7n^2-7n+2}{n^2}, \quad S^{(3)}(n) = \frac{8(n-1)^2}{n^2},$$
$$R^{(4)}(n) = 2 \frac{6n^3-9n^2+5n-1}{n^3}, \quad S^{(4)}(n) = \frac{(4n-3)(4n-4)(4n-5)}{n^3}.$$
Obviously, all coefficient functions are non-negative. 
We work out the case $r = 3$, and leave the details for $r = 4$ to the 
interested reader. By examining first few values of $x_n$, one can note
that they are slowly increasing, starting from $x_2 = 5$. Indeed, the
bounds $5 \leq x_n \leq 9$ are readily established by induction on $n$ for
$n \geq 3$. The log-convexity now follows by computing $\nabla R(n)$,
$\nabla S(n)$, $\nabla R(n-1)$, $\nabla S(n-1)$, and verifying the inequality (6)
with $m = 5$, $M = 9$.
To prove the log-balancedness of $(F^{(r)}_n)$ we start by computing
$$\Delta _{R^{(3)}}(n) = \frac{(n-1)(7n^3+7n^2-n-2)}{n^2 (n+1)^2}, \quad
\Delta _{S^{(3)}}(n) = \frac{8(n^4-2n^3-2n^2+n+1)}{n^2 (n+1)^2}.$$
It is easy to check that these determinants are positive for $n \geq 3$,
and that the conditions of Proposition 4 are valid for $n = 2$. 

Proof of the case $r = 4$ is a bit more technical, but it flows along the
same lines, and does not present any conceptual difficulties. \qed

Let us now turn our attention to the recurrences with $S(n) \leq 0$. Such
examples include, among others, Schr\"oder numbers, Delannoy numbers and, more
generally, sequences of values of Legendre polynomials. We start with a 
sequence closely connected with Franel numbers of order $3$.

The {\bf Ap\'ery numbers}, $(A_n)_{n \geq 0}$, given by the formula
$$A_n = \sum_{k = 0}^n \binom{n}{k}^2 \binom{n+k}{k}^2 =
\sum _{k = 0}^n \frac{[(n+k)!]^2}{(k!)^4 [(n-k)!]^2},$$
arose in Ap\'ery's proof of
irrationality of $\zeta (2)$ and $\zeta (3)$. They are connected with Franel
numbers of order $3$ via the identity
$$ A_n = \sum _{k = 0}^n \binom{n}{k} \binom{n+k}{k} F^{(3)}_k, \quad n \geq 0$$
(see \cite{strehl} for history of this result). The first few Ap\'ery numbers
are $1, 5, 73, 1445, 33001, 819005, \ldots $.

{\bf Corollary 4}\\ The sequence $A_n$ of Ap\'ery numbers is log-balanced.

{\bf Proof}\\ We start from the recurrence
$$A_n = \frac{34n^3 - 51n^2 +27n -5}{n^3} A_{n-1} - \frac{(n-1)^3}{n^3} A_{n-2},
\quad n \geq 2,$$
with initial conditions $A_0 = 1$, $A_1 = 5$ (\cite{beukers}). 
It is easy to prove by induction on $n$ that $x_n \geq
1$, i.e. that the sequence of Ap\'ery numbers is increasing. Hence we may take
$m = 1$ as the lower bound for $x_n$. Now the expression $\nabla R(n)
+ \nabla S(n)$ can serve as a lower bound for the expression (4), and the
log-convexity of Ap\'ery numbers follows from the inequality
$$\nabla R(n) + \nabla S(n) = \frac{1}{n^3 (n+1)^3} [50 n^4 + 52 n^3 - 10 n^2
-12 n + 4 ] \geq 0,$$
valid for all $n \geq 0$. 
For the rest, first note that $x_3 = 
\frac{1445}{73} \leq \frac{3}{2} x_2$, so we can take $n_0 = 2$. After
computing $\Delta _R (n)$ and $\overline{\Delta }_S (n)$, we get
$$\Delta _R (n) = \frac{34n^6 - 72n^4 -28n^3 +27n^2 + 7n -5}{n^3 (n+1)^3},
\quad \overline{\Delta }_S (n) = \frac{(n-1)(n^2-n-1)(2n^3+n^2-n-1)}{n^3 (n+1)^3}.$$
Both determinants are positive for $n \geq 2$, and the claim follows from
Proposition 3. \qed

{\bf Corollary 5}\\ The sequence $r_n$ of large Schr\"oder numbers is
log-balanced.

{\bf Proof}\\ Start from the recurrence
$$r_n = \frac{3(2n-1)}{n+1} r_{n-1} - \frac{n-2}{n+1} r_{n-2}, \quad n \geq 2,$$
with initial conditions $r_0 = 1$, $r_1 = 2$ \cite{stanleyII}. By computing 
the first few values of $x_n = \frac{r_n}{r_{n-1}}$, we guess the bounds
$3 \leq x_n \leq 6$, and verify them by induction for all $n \geq 2$. The
log-convexity of $(r_n)$ follows now by plugging the expressions $\nabla R(n)
= \frac{9}{(n+1)(n+2)}$ and $\nabla S(n)= - \frac{3}{(n+1)(n+2)}$ in
formula (4), together with $x_{n-1} \geq 3$. To prove the rest, we compute
$$\Delta _R(n) = 6 \frac{n-1}{n+2}, \quad \overline{\Delta }_S(n) = \frac{5-2n}{n+2}$$
and note that $\Delta _R(n) x_{n-1} + \overline{\Delta }_S(n) \geq
3 \Delta _R(n) x_{n-1} + \overline{\Delta }_S(n) \geq 0$ for all $n \geq 1$.
Hence, by Proposition 5, the sequence $(r_n)$ is log-balanced. \qed

For combinatorial interpretations of $r_n$, the reader may wish to consult Ex. 6.39 of
\cite{stanleyII}.

Our next example is the sequence of values of Legendre polynomials in some
fixed real $t \geq 1$.

{\bf Corollary 6}\\ The sequence of values of Legendre polynomials $\left (
P_n(t) \right ) _{n \geq 0} $ is log-balanced for all real $t \geq 1$.

{\bf Proof}\\ We start from Bonnet recurrence:
$$P_n(t) = \frac{2n-1}{n} t P_{n-1}(t) - \frac{n-1}{n} P_{n-2}(t), \quad n \geq 2,$$
with the initial conditions $P_0(t) = 1$, $P_1(t) = t$.
Passing to the recursion for the quotient sequence $x_n(t) = P_n(t)/P_{n-1}(t)$
we can easily establish the lower bound $x_n(t) \geq t$. By putting this 
lower bound, together with the expressions $\nabla R(n) = \frac{2}{n(n+1)}$
and $\nabla S(n)= - \frac{1}{n(n+1)}$ in formula (4), we obtain the 
log-convexity of the sequence $(P_n(t))_{n \geq 0}$. Further,
by computing $\Delta _R(n)$ and $\overline{\Delta }_S(n)$ we get
$$\Delta _R(n) = \frac{2n^2-1}{n(n+1)}, \quad \overline{\Delta }_S(n) =
\frac{-2n^2+n+1}{n(n+1)}.$$
If we suppose that $\Delta _R(n) x_{n-1}(t) + \overline{\Delta }_S(n) < 0$ for
some $n \geq 2$, we get $x_{n-1}(t) < \frac{1}{t} \frac{2n^2-n-1}{2n^2-1} <
\frac{1}{t} < t$, in contradiction with $x_{n-1}(t) \geq x_1(t) = t$.
Hence, the inequality $\Delta _R(n) x_{n-1}(t)+ \overline{\Delta }_S(n) \geq 0$
holds for all $n \geq 2$, and the claim again follows from Proposition 5. \qed

By specializing the value of $t = 3$, we get the sequence of central Delannoy
numbers, $D_n = P_n(3)$ (\cite{stanleyII}).

{\bf Corollary 7}\\ The sequence $D_n$ of central Delannoy numbers is 
log-balanced. \qed

The sequence $D_n$ counts the lattice paths from $(0,0)$ to $(n,n)$ using
only the steps $(1,0)$, $(0,1)$, and $(1,1)$. Equivalently, it counts king 
paths from the lower left to the upper right corner of an $(n+1) \times (n+1)$ 
chess board.

In all examples considered so far, the sequence $(x_n)$ was increasing, but
remained bounded. Our final example in this section shows that the same 
reasoning can be applied to the sequences whose quotient sequence increases
unboundedly.

{\bf Corollary 8}\\ Let $(a_n)$ be the sequence counting directed column-convex
polyominoes of height $n$. (See \cite{barcucci} for the definition of
these objects.) The sequence $(a_n)$ is log-balanced.

{\bf Proof}\\ From the recurrence
$$a_{n+1} = (n+1)a_n + a_1 +a_2 + \ldots + a_n, \quad n \geq 3,$$
with initial conditions $a_1 = 1$, $a_2 = 3$, given in \cite{barcucci}, one 
can easily obtain the two-term recurrence
$$a_n = (n+2) a_{n-1} - (n-1) a_{n-2}, \quad n \geq 3$$
with $a_1 = 1$, $a_2 = 3$. It can easily be shown by induction on $n$ that the
sequence $x_n = \frac{a_n}{a_{n-1}}$ is interlaced with the sequence $b_n =
n+1$, i.e. that $n+1 \leq x_n \leq n+2$. Hence the sequence $(x_n)$ is 
increasing, and $(a_n)$ is log-convex. Taking $R(n) = n+2$, $S(n) = -n + 1$,
we get $\Delta _R(n) = 2$, $\overline{\Delta }_S(n) = 1-n$. Suppose that
$\Delta _R(n) x_{n-1}(t) + \overline{\Delta }_S(n) < 0$ for some $n \geq 3$.
It follows that $x_{n-1} < \frac{n-1}{2}$, contradicting the interlacing
of $x_n$ and $b_n$. The claim now follows by checking the base of induction, i.e.
that $x_3 = \frac{13}{3} \leq \frac{3}{2} \cdot 3 = \frac{3}{2} x_2$. \qed

%We conclude the section by noting that our approach can be also successfully 
%applied to the sequences given by 
%recurrences with constant coefficients. For example, the sequence $W_n$ of
%odd-indexed Fibonacci numbers, given by $W_n = 3 W_{n-1} - W_{n-2}$, with
%$W_0 = 1$, $W_1 = 2$, is clearly log-convex, since $x_1 = 2 \leq 5/2 = x_2$,
%and the inequality $\nabla R(n) x_{n-1} + \nabla S(n) \geq 0$ is trivially valid for
%all $n \geq 2$. On the other hand, $\Delta _R(n) = 3$, $\overline{\Delta }_S(n)
%= -2$, and $\Delta _R(n) + \overline{\Delta }_S(n) \geq 0$ for all $n$, thus
%giving us the log-balancedness of $W_n$. Further, the inductive nature of
%the method becomes particularly clear here, since the relation 
%$\nabla R(n) x_{n-1} + \nabla S(n) \geq 0$ is also trivially satisfied by a
%log-concave sequence of even-indexed Fibonacci numbers, but for this sequence
%the base of induction cannot be established.

\section{Further developments}

The method exposed in Section 3 can be extended to the sequences given by a
three- (or more-) term recurrence in a straightforward way. As an illustration,
we treat here the case when all coefficient functions are positive and
increasing.

Let $(a_n)$ be a sequence of positive real numbers given by the 
recurrence $$a_n = R(n) a_{n-1} + S(n) a_{n-2} + T(n) a_{n-3}, \quad n \geq 3,$$
with given initial conditions $a_0$, $a_1$ and $a_2$. Then the recurrence for 
the quotient sequence is given by 
\be
x_n = R(n) + \frac{S(n)}{x_{n-1}} + \frac{T(n)}{x_{n-1} x_{n-2}}
\ee
for $n\geq 3$. We suppose inductively that $x_{n_0} \leq x_{n_0+1} \leq
\ldots \leq x_n$ for some $n_0 \in \N$, and we want to find sufficient 
conditions for $x_{n+1} \geq x_n$. This inequality can be stated as
$$R(n+1) + \frac{S(n+1)}{x_n} + \frac{T(n)}{x_n x_{n-1}} - R(n) - \frac{S(n)}{x_{n-1}} - \frac{T(n)}{x_{n-1} x_{n-2}} \geq 0,$$
or equivalently
$$x_n x_{n-1} x_{n-2} \nabla R(n) + x_{n-2} [x_{n-1} S(n+1) - x_n S(n)]
+ x_{n-2} T(n+1) -x_n T(n) \geq 0.$$
Now we proceed by a sequence of strengthenings of this inequality, leading to a
sufficient condition that will be expressed in known quantities and reasonably
easy to check. First we replace $S(n+1)$ and $T(n+1)$ by $S(n)$ and $T(n)$,
respectively. This yields
$$x_n x_{n-1} x_{n-2} \nabla R(n) + x_{n-2} S(n)(x_{n-1} - x_n) +T(n)
(x_{n-2}-x_{n}) \geq 0.$$
By adding $x_{n-1} - x_{n-1}$ to the term $x_{n-2}-x_{n}$ and grouping the
terms accordingly, we obtain
\be
x_n x_{n-1} x_{n-2} \nabla R(n) + [x_{n-2} S(n) +T(n)](x_{n-1} - x_n)
+T(n)(x_{n-2}-x_{n-1}) \geq 0.
\ee
Let us now look more closely at the term $x_{n-1} - x_n$. By inductive 
hypothesis, it must be non-positive, but we do not have any information about
its magnitude. Expressing $x_{n-1}$ and $x_n$ via recurrence (7) yields
$$x_{n-1} - x_n = - \nabla R(n-1) + \frac{1}{x_{n-1} x_{n-2}}
[x_{n-1} S(n-1) - x_{n-2}S(n)] + \frac{1}{x_{n-1} x_{n-2} x_{n-3}}
[x_{n-1} T(n-1) -x_{n-3} T(n)].$$
By replacing $x_{n-1}$ in the first square brackets on the right hand side of
the above relation by $x_{n-2}$, and in the second square brackets 
by $x_{n-3}$, one obtains the following inequality:
\be
x_{n-1} - x_n \geq - \nabla R(n-1) - \frac{1}{x_{n-1}}\nabla S(n-1) - \frac{1}{x_{n-1}x_{n-2}}\nabla T(n-1).
\ee
Similarly,
\be
x_{n-2} - x_{n-1} \geq - \nabla R(n-2) - \frac{1}{x_{n-2}}\nabla S(n-2) - 
\frac{1}{x_{n-2}x_{n-3}}\nabla T(n-2).
\ee
Plugging in formulae (9) and (10) in (8), we obtain the inequality
\begin{eqnarray*}
x_n x_{n-1} x_{n-2} \nabla R(n) & \geq &[x_{n-2} S(n) +T(n)]
\left [ \nabla R(n-1) + \frac{1}{x_{n-1}}\nabla S(n-1) + \frac{1}{x_{n-1}x_{n-2}}\nabla T(n-1) \right ]  \cr
  & & + T(n) \left [ \nabla R(n-2) + \frac{1}{x_{n-2}}\nabla S(n-2) + 
\frac{1}{x_{n-2}x_{n-3}}\nabla T(n-2) \right ].
\end{eqnarray*}
Finally, by replacing the values of $x_n$, $x_{n-1}$, $x_{n-2}$, and $x_{n-3}$
by their lower and upper bounds, we arrive at the following inequality:
\bea
m^3 \nabla R(n) & \geq & [M \cdot S(n) +T(n)]
\left [ \nabla R(n-1) + \frac{1}{m}\nabla S(n-1) + \frac{1}{m^2}\nabla T(n-1) \right ] \cr
& &+ T(n) \left [ \nabla R(n-2) + \frac{1}{m}\nabla S(n-2) + \frac{1}{m^2}\nabla T(n-2) \right ].
\eea
Obviously, inequality (11) implies inequality (8), and this one, in turn, 
implies our initial inequality $x_{n+1} \geq x_n$. Hence, inequality (11)
provides a sufficient condition of log-convexity for the sequence $(a_n)$.

Now, assuming the log-convexity of $(a_n)$, by following the same
reasoning as in the proof of Proposition 4, we obtain sufficient conditions
of log-balancedness of $(a_n)$ in the form
$$\Delta _R(n) x_{n-1} x_{n-2} + \Delta _S(n) x_{n-2} + \Delta _T(n) \geq 0,$$
where $\Delta _R(n)$ and $\Delta _S(n)$ are as before, and $\Delta _T(n)$ is
defined analogously.

As an illustration of this result, we prove that the sequence $(R_n)$, 
counting the 
{\bf Baxter permutations} of size $n$, is log-balanced. (See \cite{stanleyII},
p. 246 and pp. 278-9, for more details on Baxter permutations.)
The numbers $R_n$ satisfy a third-order linear recurrence with the
coefficient functions given by
$$R(n) = 2 \frac{9n^3+3n^2-4n+4}{(n+2)(n+3)(3n-2)}, \quad
S(n) = \frac{(3n-1)(n-2)(15n^2-5n-14)}{(n+1)(n+2)(n+3)(3n-2)},$$
$$T(n) = 8 \frac{(3n+1)(n-2)^2(n-3)}{(n+1)(n+2)(n+3)(3n-2)}$$
With a bit of help from a computer algebra system such as, e.g. {\it
Mathematica}, it can be proved that $7 \leq x_n \leq 9$ for $n \geq 47$. 
Verifying the inequality (11) then boils down to checking that a certain 
rational function of $n$ (with the degrees of the numerator and denominator
equal to $12$ and $14$, respectively) is nonnegative for sufficiently large
values of the argument. By substituting $n+3$ in place of $n$ it becomes
obvious that all the coefficients become positive, and hence, the function
cannot change the sign for $n \geq 3$. The increasing behavior of $x_n$ for 
$n \leq 47$ is
easily checked by direct computation. Hence the sequence $(R_n)$ is log-convex.
To prove the log-balancedness, it is easy to check that all three determinants
$$\Delta _R(n) = \frac{27n^5+18n^4+3n^3+76n^2+100n+16}{(n+1)(n+2)(n+3)(n+4)(3n+1)(3n-2)},$$
$$\Delta _S(n) = \frac{135n^5-990n^4+87n^3+1036n^2+4n-112}{(n+1)(n+2)(n+3)(n+4)(3n+1
)(3n-2)},$$
$$\Delta _T(n) = \frac{9n^5-138n^4+349n^3-80n^2-252n-48}{(n+1)(n+2)(n+3)(n+4)(
3n+1
)(3n-2)}$$
are positive for $n \geq
13$, and the log-balancedness of $(R_n)$ follows by directly verifying 
defining inequalities in the remaining cases.
All the {\it Mathematica} calculations necessary for verifying the above
inequalities were performed exactly.

The scope of our approach can also be extended in another direction, namely
to linear nonhomogeneous recurrences. Here we indicate, after the fashion 
of \cite{dv}, how such recursions
can be transformed in a form suitable for application of our method. So, for 
example, let $(a_n)$ be given by a linear nonhomogeneous recurrence of the
first order
\be
a_n = R(n) a_{n-1} + S(n)
\ee
with the initial condition $a_0$.
By writing down the recurrence (12) for successive indices, multiplying and
subtracting as to cancel the nonhomogeneous part, one obtains the homogeneous
second order linear recurrence for $a_n$:
%Consider the quotients $q(n)=\frac{a(n)}{a(n-1)}$ and note that
%\be
%a(n)=q(n)q(n-1) \ldots q(2)a(1), \quad n\geq 2.
%\ee
%Then, dividing (12) by $a(n-1)$ we obtain a (long) recurrence for $q(n)$'s:
%\be
%q(n)=R(n)+\frac{S(n)}{q(n-1)q(n-2) \ldots q(2)a(1)}.
%\ee
%To get a short recurrence for $q(n)$'s, substitute for $a(n)$ and $a(n-1)$ the
%corresponding products (13) in (12) ($n\geq 3$):
%$$\displaylines{
%\qquad q(n)q(n-1)\ldots q(2)a(1)=R(n)q(n-1)\ldots q(2)a(1)+S(n)= \hfill \cr
%\hfill R(n)\frac{q(n)q(n-1)\ldots q(2)a(1)}{q(n)}+S(n)\frac{q(n)q(n-1)\ldots q(2)a(1)}{q(n)q(n-1)\ldots q(2)a(1)}. \qquad\cr
%}$$
%From there we get
%$$\frac{1}{q(n)q(n-1)\ldots q(2)a(1)}=\frac{1}{S(n)}\left [ 1-\frac{R(n)}{q(n)} \right ],$$
%and then
%\be
%\frac{1}{q(n-1)\ldots q(2)a(1)}=\frac{1}{S(n-1)}\left [ 1-\frac{R(n-1)}{q(n-1)} \right ] .
%\ee
%Substituting (15) in (14) yields a short recursion for $q(n)$'s:
$$a_n=\left [ R(n)+\frac{S(n)}{S(n-1)} \right ] a_{n-1}-\frac{R(n-1)S(n)}{S(n-1)}a_{n-2}.$$
By denoting $R^{\star}(n) = R(n)+\frac{S(n)}{S(n-1)}$, 
$S^{\star}(n) = -\frac{R(n-1)S(n)}{S(n-1)}$, and dividing through by $a_{n-1}$,
we get a recurrence for $x_n$ of the type (2) and the further treatment 
depends on the combination of signs of $R^{\star}(n)$ and $S^{\star}(n)$.

Similarly, for a second order linear recurrence
$$a_n=R(n)a_{n-1}+S(n)a_{n-2}+T(n),$$
we obtain
$$x_n=R(n)+\frac{S(n)}{x_{n-1}}+\frac{T(n)}{T(n-1)}\left [ 1- \frac{R(n-1)}{x_{n-1}}-\frac{S(n-1)}{x_{n-1}x_{n-2}} \right ]. $$
Then we can proceed as before.

Finally, a word of caution. It would be hasty to conclude, from the 
cited examples, that all combinatorially interesting sequences are log-balanced.
For example, the sequences $a_n = (n!)^2$, $a_n = (n-1)!$ and 
$a_n = \sum _{k = 0}^n k!$
are not log-balanced, since their quotient sequences grow too fast. It is also
interesting to note that the property of log-balancedness is not shift-invariant;
one can easily see that the sequence $(n+1)!$ is log-balanced, while $(n-1)!$
is not.

One could, in principle, consider an alternative approach to the question of
log-balancedness, that is in a sense dual to ours. One could take a 
log-concave sequence $(a_n)$ and ask for the sufficient conditions for the
sequence $(n!a_n)$ to be log-convex. Since it appears that the log-convex
sequences are much more common among the sequences of combinatorial interest,
we will not pursue this alternative approach here.

The author acknowledges the support of the Welch Foundation of Houston, Texas,
via grant \# BD-0894.

\end{document}